\documentclass[12pt]{article}

\pdfoutput=1
\usepackage{amsfonts}
\usepackage{amsmath}
\usepackage{amssymb}
\usepackage{amsthm}
\usepackage{graphicx}
\usepackage{color}
\usepackage{subfigure}
\usepackage[hcentering,vcentering]{geometry}  
\newcommand{\field}[1]{\mathbb{#1}}
\newcommand{\R}{\field{R}}
\newcommand{\N}{\field{N}}
\newcommand{\bu}{{\bf u}}

\newcommand{\bx}{{\bf x}}
\newcommand{\bw}{{\bf w}}

\newcommand{\cP}{{\cal P}}
\newcommand{\cF}{{\cal F}}

\newcommand{\cQ}{{\cal Q}}
\newcommand{\cB}{{\cal B}}
\newcommand{\cX}{{\cal X}}
\newcommand{\cG}{{\cal G}}
\newcommand{\cN}{{\cal N}}
\newcommand{\cD}{{\cal D}}

\newcommand{\cO}{{\cal O}}
\newcommand{\eps}{{\varepsilon}}

\newtheorem{corollary}{Corollary}
\newtheorem{definition}{Definition}
\newtheorem{algorithm}{Algorithm}
\newtheorem{proposition}{Proposition}
\newtheorem{theorem}{Theorem}
\newtheorem{remark}{Remark}

\def\KK{\mathcal{K}}
\def\LL{\mathcal{L}}
\def\argmin{\mathop{\rm argmin}\nolimits}

\def\endproof{\unskip \nobreak \hskip0pt plus 1fill \qquad \qed}

\begin{document}

\title{Global optimal control of perturbed systems}

\author{Lars Gr\"une\footnote{Mathematisches Institut, 
Universit\"at Bayreuth,
95440 Bayreuth, Germany, 
{\tt lars.gruene@uni-bayreuth.de}}
\and
Oliver Junge\footnote{
Zentrum Mathematik,
Technische Universit\"at M\"unchen,
85748 Garching, Germany,
{\tt junge@ma.tum.de}}}

\date{October 2006}

\maketitle

\begin{abstract}
We propose a new numerical method for the computation of the optimal value function of perturbed control systems and associated globally stabilizing optimal feedback controllers.  The method is based on a set oriented discretization of state space in combination with a new  algorithm for the computation of shortest paths in weighted directed hypergraphs.  Using the concept of a multivalued game, we prove convergence of the scheme as the discretization parameter goes to zero.

\medskip

\end{abstract}

\section{Introduction}

Global infinite horizon optimal control methods for the solution of
general nonlinear stabilization problems are attractive for their
flexibility and theoretical properties, because they are applicable to
virtually all types of nonlinear dynamics, their optimal value functions
can typically be identified as Lyapunov functions and they allow for a
rigorous treatment of perturbations in a game theoretical
setting. However, these methods have the drawback that their numerical
solution requires the discretization of the state space which results
in huge numerical problems both in terms of computational cost and in
terms of memory requirements. Hence, in order to make these methods
applicable to a broader range of systems, advanced numerical
techniques are needed in order to reduce the computational effort as
much as possible. 

A novel approach to such problems was presented in the recent paper
\cite{JungeOsinga:04}, where a set oriented numerical method for the
approximate computation of the optimal value function of certain
nonlinear optimal control problems has been developed.  The approach
relies on a division of state space into boxes that constitute the
nodes of a directed weighted graph, where the weights are constructed
from the given cost function.  On this graph, standard graph theoretic
algorithms for computing shortest paths can directly be applied,
yielding an approximate value function which is piecewise constant on
the state space.  At the same time, for every node in the graph, these
algorithms compute the successor node on a shortest path, yielding
approximate optimal pseudo-trajectories of the original system.
Hence, this method combines a simple and hierarchically implementable 
discretization technique with efficient graph theoretic algorithms
yielding both low memory consumption and a fast solution.
For the problem of feedback stabilization the solution from
\cite{JungeOsinga:04}, however, is not directly applicable, because
the resulting pseudo-trajectories would have to be postprocessed in
order to obtain true solutions of the system.    

In \cite{GrJu:05} it was subsequently shown that the approximate
optimal value function can in fact be used in order to construct a
stabilizing feedback controller.  Based on concepts from dynamic
programming \cite{Gruene:97} and Lyapunov based approximate stability
analysis \cite{NesicTeel:04},  a statement about its optimality
properties was given and a local a posteriori error estimate derived
that enables an adaptive construction of the division of state space.
However, due to the fact that the approximate optimal value function
is not continuous, the constructed feedback law is in general not
robust with respect to perturbations of the system. 

In the present paper, we show how to incorporate arbitrary
perturbations into the framework sketched above. These perturbations
can be either inherently contained in the underlying model,
describing, e.g., external disturbances or the effect of unmodelled
dynamics, or they could be added on top of the original model to
account, e.g., for discretization errors.

Our goal in this paper is to
construct a feedback which is robust in the sense that on a certain
subset of state space it stabilizes the system regardless on how the
perturbation acts.  Conceptually, this problem leads to a dynamic
game, where the controls and the perturbations are associated to two
``players'' that try to minimize and to maximize a given cost
functional, respectively.  We show how the discretization of state
space in a natural way leads to a multivalued dynamic game (i.e.\ a
discrete inclusion) and prove convergence of the associated value
function when the images of the inclusion shrink to the original
single-valued map.  From this multivalued game we derive a directed
weighted hypergraph that gives a finite state model of the original
game.   We formulate an adapted version of Dijsktra's algorithm in
order to compute the associated approximate value function and prove
convergence when the box-diameter of the state space division goes to
zero. 

It should be noted that the convergence analysis developed in this
paper using multivalued dynamics is new also for the discretization of
optimal control problems without perturbations in
\cite{JungeOsinga:04}. An interesting side result of our study is
that using this technique we are able to keep track of the effects of 
discontinuities in the approximated optimal value function as induced,
e.g., by state space constraints. This allows us to prove not only
$L^\infty$ convergence in regions of continuity but also $L^1$
convergence in the whole domain of the optimal value function,
provided that the optimal value function is continuous with respect to
small changes in the state space constraints. 

Compared to other dynamic programming approaches to the stabilization
of perturbed nonlinear systems (see, e.g., \cite{HuJaNeDo:05} and the
references therein), the main advantages of our method are these
general and rigorously provable convergence properties and the low
computational cost of our perturbed version of Dijkstra's algorithm,
cf.\ Section~\ref{ssec:implementation}. However, 
our new algorithm is also advantageous for unperturbed problems when
treating the spatial discretization errors as perturbation: as Example
\eqref{eq:pendulum} illustrates, this approach leads to considerably
improved performance on a significantly coarser discretization
compared to \cite{GrJu:05}. 

The paper is organized as follows.  In the ensuing
Section~\ref{sec:problem} we describe the problem formulation and the
associated game theoretic interpretation.  In Section~\ref{sec:mvg} we
introduce the concept of a multivalued game and an enclosure and prove
a statement about the convergence of the value function of a sequence
of enclosures of a multivalued game. These result are extended to
systems with state constraints in Section~\ref{stateconstraints}. In
Section~\ref{sec:discretization} we show how via the division of state
space one obtains  a multivalued game from the original system,
construct the corresponding hypergraph and introduce an associated
shortest path algorithm. Some hints on its implementation, complexity
issues as well as two numerical examples are addresed in 
Section~\ref{sec:implementation}. Convergence of the numerical
approximation to the optimal value function and the construction of
approximately optimal feedback laws are discussed in Sections
\ref{sec:convergence} and \ref{sec:feedback}, respectively.

\section{Problem formulation}
\label{sec:problem}

We consider the problem of optimally stabilizing the discrete-time 
perturbed control system
\begin{equation}
  \label{eq:evol}
  x_{k+1} = f(x_k, u_k, w_{k}), \quad k=0,1,\ldots,
\end{equation}
where $f:X\times U\times W\to X$ is continuous, $x_{k}\in X$ is the
state of the system, $u_{k}\in U$ is the control input and $w_{k}\in
W$ is a perturbation parameter,  chosen from sets $X\subset\R^{d},
U\subset\R^m$ and $W\subset\R^\ell$. In addition to the evolution law,
we are given a continuous cost function $g:X\times U\to [0,\infty)$,
that assigns the cost $g(x_{k},u_{k})$ to any transition 
$x_{k+1} = f(x_k, u_k, w_{k})$, $w_{k}\in W$. 

Our goal is to derive an (optimal) \emph{feedback law} $u:X\to U$
that \emph{stabilizes} the system in the sense that for a certain
subset $S\subset X$ any trajectory starting in $S$ tends to some
prescribed set $O\subset X$, while the \emph{worst case
accumulated cost} is minimized. 

Let us be more precise.  For a given initial point $x\in X$, a control sequence $\bu=(u_{k})_{k\in\N}\in U^\N$ and a perturbation sequence
$\bw=(w_{k})_{k\in\N}\in W^\N$ yield the \emph{trajectory}
$\bx(x,\bu,\bw)=(x_{k}(x,\bu,\bw))_{k\in\N}$, defined by 
$x_{0}=x$ and
\begin{equation}
\label{pcs}
x_{k+1}=f(x_{k}(x,\bu,\bw), u_{k}, w_{k}), \quad k=0,1,\ldots, 
\end{equation}
while the associated accumulated cost is given by
\[
J(x,\bu,\bw)=\sum_{k=0}^\infty g(x_{k}(x,\bu,\bw),u_{k}).
\]

In order to formalize the interplay between the control and the
perturbation we employ a \emph{game theoretic viewpoint} which we
describe next. 
The problem formulation actually already describes a game (see,
e.g., \cite{F:61}), where at each step of the iteration
(\ref{eq:evol}) two ``players'' choose a control value $u_{k}$ and a
perturbation value $w_{k}$, respectively.  The goal of the controlling
player is to minimize $J$, while the perturbing player tries to
maximize this quantity. 

We assume that the controlling player has to choose the value $u_{k}$
first and that the perturbing player has the advantage of  knowing
$u_{k}$ when choosing the perturbation value $w_{k}$.   However, the
perturbing player is not able to forsee future choices of the
controlling one.  More formally, we restrict the choice of
perturbation sequences $\bw\in W^\N$ to those that result from
applying a \emph{nonanticipating strategy} $\beta:U^\N\to W^\N$ to a
given control sequence $\bu\in U^\N$, i.e. we have $\bw=\beta(\bu)$,
with $\beta$ satisfying 
\[
u_{k}=u_{k}' \quad\forall k\le K\quad\Rightarrow\quad
\beta(\bu)_{k}=\beta(\bu')_{k} \quad\forall k\le K 
\]
for any two control sequences 
$\bu=(u_{k})_{k},\bu'=(u'_{k})_{k}\in U^\N$.  Let $\mathcal{B}$ denote
the set of all nonanticipating strategies $\beta:U^\N\to W^\N$. 

As mentioned, our goal is to find a feedback law $u:X\to U$ such
that with controls $u_{k}=u(x_{k})$, $x_{k}$ approaches a given set $O\subset X$, regardless of how the perturbation sequence $\bw$ is
chosen. 
% Depending on the way that the perturbation parameter $w_{k}$
%acts in $f$, it may be impossible that $x_{k}$ actually tends to $0$
%as $k\to\infty$.  However, in an application context it may be
%sufficient to steer $x_{k}$ close to the origin and this is the
%situation we have in mind. 
Accordingly, we assume that we know a compact robust
forward invariant set $O\subset X$, i.e.\
for all $x\in O$ there is a control $u\in U$ such that
$f(x,u,W)\subset O$.  Since we are done with controlling the
system once we are on $O$, we assume that $g(x,u)=0$ for all
$x\in O$ and all $u\in U$ and $g(x,u)>0$ for all
$x\not\in O$ and all $u\in U$. Further assumptions on $g$ and on the
dynamics in a neighborhood of $O$ will be specified later.  

Our construction of the feedback law will be based on the \emph{upper
value function} $V:X\to[0,\infty]$, 
\begin{equation}
\label{valuefunction}
V(x)=\sup_{\beta\in\mathcal{B}}\inf_{\bu\in U^\N} J(x,\bu,\beta(\bu)),
\end{equation}
of the game~(\ref{eq:evol}), which fulfills the \emph{optimality principle}
\begin{eqnarray}
\label{eq:optimality principle}
V(x) &=& \inf_{u\in U} \left[ g(x,u) + \sup_{w\in W} V(f(x,u,w))\right].
\end{eqnarray}

\section{Multivalued games}
\label{sec:mvg}

As we will see in the next section, our set oriented approach to
the discretization of state space of the perturbed control
system~(\ref{eq:evol}) leads to a finite state multivalued system. 
For the convergence analysis of this discretization it turns out to be
useful to introduce as an intermediate object an infinite state
multivalued game defined by a discrete inclusion. This is given by a
multivalued map
\[
F:X\times U\times W\rightrightarrows X,
\]
where $X\subset\R^d$ is a closed set and $U\subset\R^m$, $W\in\R^\ell$
and the images of $F$ are compact sets, together with a cost function
\[
G:X\times X\times U \times W\to [0,\infty).
\]
In order to simplify our presentation we first assume that 
$F(x,u,w)\ne \emptyset$ for all $x\in X$, $u\in U$, $w\in W$, which
will be relaxed later, cf.\ Section~\ref{stateconstraints}. Further regularity assumptions on these maps will be imposed when 
needed. Note that we have introduced a second state argument in $G$,
which allows to associate different costs to the trajectories of the
associated discrete inclusion. 

For a given initial state $x\in X$, a given control sequence
$\bu=(u_{k})_{k\in\N}\in U^\N$ and a given perturbation sequence
$\bw=(w_{k})_{k\in\N}\in W^\N$, a \emph{trajectory} of the game is
given by any sequence $\bx=(x_{k})_{k\in\N}\in X^\N$ such that
$x_{0}=x$ and 
\[
x_{k+1}\in F(x_{k},u_{k}, w_{k}), \quad k=0,1,2,\ldots.
\]
We denote by
\[
\mathcal{X}_{F}(x,\bu,\bw)=\left\{(x_{k})_{k}\in X^\N\mid x_{0}=x, x_{k+1}\in F(x_{k},u_{k},w_{k}) \;\forall k\in\N\right\}
\]
the set of all trajectories of $F$ associated to $x$, $\bu$ and $\bw$.
The \emph{accumulated cost} is given by
\[
J_{(F,G)}(x,\bu,\bw) = \inf_{(x_{k})_{k}\in\mathcal{X}_{F}(x,\bu,\bw)}\sum_{k=0}^\infty G(x_{k},x_{k+1},u_{k},w_{k}).
\]
As in the previous section, we are interested in computing the upper
value function   
\begin{equation}
\label{multival}
V_{(F,G)}(x)=\sup_{\beta\in\cB}\inf_{\bu\in U^\N} 
J_{(F,G)}(x,\bu,\beta(\bu)), \quad x\in X,
\end{equation}
of this game. By standard dynamic programming arguments
\cite{Bertsekas:95} one sees that this function fulfills the
optimality principle 
\begin{equation}
\label{multi-optprinc}
V_{(F,G)}(x)=\inf_{u\in U}\sup_{w\in W}\inf_{x_{1}\in F(x,u,w)}\left\{
G(x,x_{1},u,w) + V_{(F,G)}(x_{1})\right\}. 
\end{equation}

Observe that our original ``single valued'' game
(\ref{pcs})--(\ref{valuefunction}) can be
recast in this multivalued setting by defining
\[ F(x,u,w) := \{f(x,u,w)\} \mbox{ and } G(x,x_1,u,w) := g(x,u). \]

We will now investigate the relation of the value functions of
different multivalued games. For this purpose we first introduce the
concept of an \emph{enclosure}.

\begin{definition}
If $(F_{1},G_{1})$ and $(F_{2},G_{2})$ are two multivalued games such that 
\[
F_{2}(x,u,w) \subset F_{1}(x,u,w) 
\]
for all $x, u$ and $w$ and 
\[
G_{1}(x,x',u,w) \leq G_{2}(x,x',u,w)
\]
for all $x, x'\in F_{2}(x,u,w)$ and all $u$ and $w$, then
$(F_{1},G_{1})$ is called an \emph{enclosure} of $(F_{2},G_{2})$. 
\end{definition}

From this definition we immediately obtain the following proposition.

\begin{proposition}\label{prop:enclosure}
Let the game $(F_{1},G_{1})$ be an enclosure of the game
$(F_{2},G_{2})$.  Then  
\[
V_{(F_{1},G_{1})}\leq V_{(F_{2},G_{2})}.
\]
\end{proposition}

The next proposition studies the convergence of the value functions
$V_{(F_i,G_i)}$ of a sequence of games $(F_i,G_i)$. In this
proposition $H$ denotes the Hausdorff distance for compact sets.

\begin{proposition}\label{prop:convergence} Let the sequence of games 
$(F_{i},G_{i})$, $i\in\N$, be enclosures of the game $(F,G)$ and assume
\begin{equation} \sup_{x\in X,u\in U, w\in W} H(F_i(x,u,w),F(x,u,w))
\to 0 \quad \mbox{as } \,  i\to \infty
\label{Fconvergence}\end{equation} 
and 
\begin{equation} \sup_{x,x_1\in X,u\in U, w\in W} |G_{i}(x,x_1,u,w) -
G(x,x_1,u,w)| 
\to 0 \quad \mbox{as } \,  i\to
\infty. \label{Gconvergence}\end{equation} 
Assume furthermore that $F$ is upper semi--continuous in $x$ and that
$G$ is continuous in $x$ and $x_1$, both uniformly 
in $u$ and $w$ and on compact subsets of $X$. 
In addition, we assume that there exists $\alpha\in\KK_\infty$\footnote{A function
  $\gamma: [0,\infty)\to [0,\infty)$ is of class
  $\KK$ if it is continuous, zero at zero and strictly increasing.
  It is of class $\KK_{\infty}$, if, in addition, it is unbounded.} with
\[ G(x,x_1,u,w) \ge \alpha(d(x,O)+d(x_1,O)) \]
and 
\[ G_i(x,x_1,u,w) \ge \alpha(d(x,O) + d(x_1,O)) \]
for all $i\in\N$, $u\in U$, $w\in W$, and that $V_{(F,G)}$ is
continuous on $\partial O$. Then for each compact set $K\subset X$
for which $\sup_{x\in K} V_{(F,G)}(x) < \infty$ we have  
\[ \sup_{x\in K} |V_{(F_i,G_i)}(x) -
V_{(F,G)}(x)| \to 0 \quad 
\mbox{as } \,  i\to \infty, \]
i.e., uniform convergence on compact sets in the domain of 
$V_{(F,G)}$.
\end{proposition}

\begin{proof}
Let $k^*:X^\N \to \N$ be a bounded map. Then 
from the optimality principle (\ref{multi-optprinc}) we obtain by
induction 
\begin{eqnarray*} && V_{(F,G)}(x) =
\sup_{\beta\in\mathcal{B}}\inf_{\bu\in U^\N}
  \inf_{\bx\in\mathcal{X}_{F}(x,\bu,\beta(\bu))}
\Bigg\{\sum_{k=0}^{k^*(\bx)-1} 
  G(x_{k},x_{k+1},u_{k},\beta(\bu)_{k}) \\ 
&& \hspace*{8cm} + \enspace V_{(F,G)}(x_{k^*(\bx)})\Bigg\}
\end{eqnarray*}

Now let $\gamma:=\sup_{x\in K} V_{(F,G)}(x)$. 
Due to the lower bound $\alpha$ on $G$, for every $\delta >0$  there exists a time $k_{\gamma,\delta}\in\N$ such that for  each trajectory $\bx\in\cX_F(x,\bu,\beta(\bu))$ with cost bounded by
$\gamma$ there exists a time
$k^*(\bx) \le k_{\gamma,\delta}$ such that 
$x_{k^*(\bx)}\in B_\delta(O)$.
We fix $\eps>0$ and $x\in K$ and choose
$\delta>0$ such that $V_{(F,G)}(x)\le \eps$ for all 
$x\in B_\delta(O)$ ($\delta$ exists because of the continuity of $V_{(F,G)}$ on
$\partial O$). Then, using an $\eps$--optimal perturbation strategy
$\beta^*\in\cB$ and an arbitrary $\bu^*\in U^\N$, from the above optimality
principle we obtain  
\begin{eqnarray*}
 V_{(F,G)}(x) 
& \le & \inf_{\bu\in U^\N}
  \inf_{\bx\in\mathcal{X}_{F}(x,\bu,\beta^*(\bu))}
  \left\{\sum_{k=0}^{k^*(\bx)-1} 
  G(x_{k},x_{k+1},u_{k},\beta^*(\bu)_{k}) \right.\\
&& \hspace*{6cm} + V_{(F,G)}(x_{k^*(\bx)})\Bigg\}
    + \eps \\ 
& \le & \inf_{\bu\in U^\N}
  \inf_{\bx\in\mathcal{X}_{F}(x,\bu,\beta^*(\bu))}
  \left\{\sum_{k=0}^{k^*(\bx)-1} 
  G(x_{k},x_{k+1},u_{k},\beta^*(\bu)_{k})\right\} + 2\eps \\ 
& \le & \inf_{\bx\in\mathcal{X}_{F}(x,\bu^*,\beta^*(\bu^*))}
\left\{\sum_{k=0}^{k^*(\bx)-1} 
  G(x_{k},x_{k+1},u^*_{k},\beta^*(\bu^*)_{k})\right\} + 2\eps.
\end{eqnarray*}

Now, fixing $\beta^*$, for any $i\in\N$ we can pick an
$\eps$--optimal control $\bu_i^*$, yielding 
\begin{eqnarray*} \gamma & \ge & V_{(F_i,G_i)}(x) \\
& \ge & \inf_{\bx\in\mathcal{X}_{F_i}(x,\bu_{i}^*,\beta^*(\bu_i^*))}
  \left\{\sum_{k=0}^{\infty}
  G_i(x_{k},x_{k+1},(\bu_i^*)_{k},\beta^*(\bu_i^*)_{k})\right\}- \eps \\
& \ge & \inf_{\bx\in\mathcal{X}_{F_i}(x,\bu_{i}^*,\beta^*(\bu_i^*))}
  \left\{\sum_{k=0}^{k^*(\bx)}
  G_i(x_{k},x_{k+1},(\bu_i^*)_{k},\beta^*(\bu_i^*)_{k})\right\} - \eps.
\end{eqnarray*}
In particular, this last expression is bounded by $\gamma$ and hence 
the lower bound $\alpha$ for $G_i$ 
implies that there exists a compact set $K_1$ such
that each $\eps$--optimal trajectory
$(x_k)_k\in\mathcal{X}_{F_i}(x,\bu_{i}^*,\beta^*(\bu_i^*))$ lies in $K_1$ for
all $i\in\N$. 

Now assumption (\ref{Fconvergence}) and the upper semicontinuity of
$F$ imply that for each $\eps_1>0$ 
there exists an $i_0\in\N$ such that for $i\geq i_{0}$ and each such $\eps$--optimal
trajectory $(x_k)_k\in\mathcal{X}_{F_i}(x,\bu_{i}^*,\beta^*(\bu_i^*))$ there
exists a trajectory $(\tilde x_k)_k \in \cX_{F}(x,\bu_{i}^*,\beta^*(\bu_{i}^*))$ with
$\|x_k-\tilde x_k\|\le \eps_1$ for all
$k=1,\ldots,k_{\gamma,\delta}$. Hence (\ref{Gconvergence}) and the
continuity of $G$ imply that we can
find $i_1\in\N$ such that  
\begin{eqnarray*}
&& \left|\inf_{(x_{k})_{k}\in\mathcal{X}_{F}(x,\bu_{i}^*,\beta^*(\bu_i^*))}
  \left\{\sum_{k=0}^{k^*}
  G(x_{k},x_{k+1},(\bu_i^*)_{k},\beta^*(\bu_i^*)_{k})\right\}\right. \\
   && \left.\quad -
  \inf_{(x_{k})_{k}\in\mathcal{X}_{F_i}(x,\bu_{i}^*,\beta^*(\bu_i^*))} 
  \left\{\sum_{k=0}^{k^*}
  G_i(x_{k},x_{k+1},(\bu_i^*)_{k},\beta^*(\bu_i^*)_{k})\right\} 
  \right| \le \eps
  \end{eqnarray*} 
for all $i\ge i_1$ and all
$k^*\in\{1,\ldots,k_{\gamma,\delta}\}$. Combining this inequality with
the estimates for 
$V_{(F,G)}$ and $V_{(F_i,G_i)}$ using $\bu^*=\bu_i^*$ in the former we
obtain 
\[ V_{(F,G)}(x) \le V_{(F_i,G_i)}(x) + 5\eps \]
for all $i\ge i_1$. Since $i_1$ depends only on $k_{\gamma,\delta}$
and $\eps$, hence only on the set $K$ and not on the individual $x$,
we obtain the desired uniform convergence. 
\end{proof}

\begin{remark} Note that we have obtained our result under very weak
assumptions on $F$ and $G$ using, however, the crucial continuity
assumption of $V_{(F,G)}$ on $\partial O$. This assumption --- which is
implicit and in general difficult to check directly --- can be ensured
by the following asymptotic controllability assumption on the dynamics
$F$ and the cost function $G$ in a neighborhood of $O$: 

Assume that there exists a neighborhood $\mathcal{N}$ of $O$ and a
$\KK\LL$ function\footnote{A function $\eta:[0,\infty)\times[0,\infty)\to[0,\infty)$ is of
  class $\KK\LL$ if it is continuous, of class $\KK$ 
  in the first variable and strictly decreasing to $0$ in the second
  variable.} $\eta$ such that for each $x\in\mathcal{N}$ and
each perturbation strategy $\beta\in\cB$ there exists a
control sequence $\bu\in U^\N$ 
and a trajectory $(x_{k})_{k}\in\mathcal{X}_{F}(x,\bu,\beta(\bu))$
with 
\begin{equation} d(x_k,O)\le \eta(d(x_0,O),
k). \label{eq:ascontrol}\end{equation} 
Then, using the construction from \cite[Proof of Theorem
5.4]{GrueneNesic:03}, we find a $\KK$ function $\rho$ (denoted 
$\rho_2$ in \cite{GrueneNesic:03}) such that 
$G(x_{0},x_1,u,w)\le \rho(d(x_{0},O))$ for $x_{0}\in\mathcal{N}$ implies 
\[ \sum_{k=0}^\infty G(x_{k},x_{k+1},u_{k},\beta(\bu)_{k}) \le
\tilde\sigma(d(x_{0},O)) \]
for some $\KK$ function $\tilde\sigma$. Since $\tilde\sigma(d(x,O)) \to 0$ as $d(x, O)\to 0$ this implies $V(x)\to 0$ as $d(x, O)\to 0$ which yields
continuity of $V$ on $\partial O$. Note that condition
\eqref{eq:ascontrol} is weaker
than controllability conditions typically employed to ensure
continuity in minimum time problems or pursuit--evasion games (cf.\
e.g.\ \cite[Chapter IV]{BardiCapuzzoDolcetta97}) because we do not
require to be able to steer the system {\em into} the ``target'' set
$O$ but only {\em asymptotically to} $O$.

We also emphasize that we only need continuity at the boundary of $O$
and that our optimal value function may be discontinuous elsewhere.
\end{remark}

\section{State space constraints}\label{stateconstraints}

So far we have assumed $F(x,u,w)\ne \emptyset$ for all $x\in X$, $u\in
U$, $w\in W$ which guarantees that for each initial value $x$, and
each pair of control and perturbation sequences $\bu$ and $\bw$ we
obtain at least one trajectory $(x_k)_k$ which is defined for all
$k\in\N_0$. However, in practice it will often be necessary to relax
this assumption. 

In order to motivate this relaxation, assume that we are given a
multivalued game $(\widetilde F, G)$ on a state space 
$\widetilde X\subseteq \R^d$. In our numerical approach, the state
space set $X$ on which we can solve the problem will be a compact set
while the state space $\widetilde X$ of the given 
problem is often unbounded. In addition, from a modeling point of
view it might be desirable to introduce state constraints,
e.g., in order to avoid certain critical regions of the state space. 
In both cases, it will be necessary to restrict the state space of
the original problem defining
\[ F(x,u,w) := \widetilde F(x,u,w)\cap X, \enspace 
x\in X,\, u\in U,\, w\in W. \]
This construction may result in $F(x,u,w) = \emptyset$ for
certain 
$x\in X$, $u\in U$, $w\in W$ and consequently it may happen that a
solution trajectory will only exist for finite time.  More precisely,
for given $F$, given $\bu=(u_{k})_{k}\in U^\N$, given
$\bw=(w_{k})_{k}\in W^\N$ and any sequence $\bx=(x_{k})_{k}\in X^\N$
let  
\[
k_F^{\max}(\bx,\bu,\bw)=\max\;\left\{ \hat k\in\N: x_{k+1}\in
F(x_{k},u_{k},w_{k}), k=0,\ldots,\hat k-1\right\} 
\]
be the maximal index up to which the sequence $\bx$ constitutes a
trajectory of $F$.  Since a trajectory with $k_{F}^{\max}(\bx,\bu,\bw)
< \infty$ cannot converge to the set $O$ we set 
\[ 
J_{(F,G)}(x,\bu,\bw) := \infty \quad \mbox{ if } 
k_F^{\max}(\bx,\bu,\bw)<\infty \mbox{ for each }
\bx\in X^\N \mbox{ with } x=x_{0}.
\]
It is easy to see that
Proposition \ref{prop:enclosure} remains valid in this case, while
Proposition \ref{prop:convergence} is more 
difficult to recover in this setting. The reason lies in the fact that
any enclosure will necessarily enlarge the set of possible 
trajectories, even if we apply the same state space constraints to $F$
and $F_i$. In the presence of state space constraints this means that
for any $i$ there may exist a trajectory $\bx_i$ of $F_i$ for which all
nearby trajectories $\bx$ of $F$ violate the space constraints. 
In other words, unless very specific knowledge about the dynamics $F$
is available and used for the construction of the enclosure $F_i$, the
enlargement of the dynamics has the implicit effect of relaxing the
state space constraints.  

However, if we assume that the optimal value function is continuous
with respect to relaxations of the state space
constraints, then we can recover Proposition
\ref{prop:convergence}. In order to formalize this relaxation, for
$\eps>0$ we define the space 
\[ X_\eps:=\{ x\in \widetilde X\,|\, d(x,X)\le \eps \}, \]
the multivalued dynamics 
\[ F_\eps(x,u,w):=\widetilde F(x,u,w)\cap X_\eps \]
and the related optimal value function $V_{(F_\eps,G)}$. 
Using this notation we can prove the following variant of Proposition
\ref{prop:convergence}.

\begin{proposition}\label{prop:sc-convergence} Consider the state
space constrained dynamics $F$ of $\tilde F$ and consider a sequence of
enclosures $(F_i,G_i)$ of $F$ on $X$.
Let the assumptions of Proposition \ref{prop:convergence} hold for
$F$ and $F_i$, where (\ref{Fconvergence}) in the case of
$F(x,u,w)=\emptyset$ is to be understood as 
\[
 F_i(x,u,w)=\emptyset \mbox{ for all }
i\in\N \mbox{ and all } x,u,w \mbox{ with }
F(x,u,w)=\emptyset.
\]
Assume, furthermore, that $\widetilde F$ is upper semi--continuous in $x$
uniformly in $u$ and $w$ on compact subsets of $\widetilde X$ and let
$\|\cdot\|_p$ be the usual $p$--norm for real valued functions on $X$
for some $p\in\{1,\ldots,\infty\}$. 

Then for each compact set $K\subset X$
for which $\sup_{x\in K} V_{(F,G)}(x) < \infty$ and on which the
continuity assumption
\begin{equation} \|V_{(F_\eps,G)}|_K-V_{(F,G)}|_K\|_p \to 0 \enspace \mbox{as}
\enspace \eps \to 0 \label{Vcont}\end{equation}
holds, we have  
\[ \|V_{(F_i,G_i)}|_K - V_{(F,G)}|_K\|_p \to 0 \quad 
\mbox{as } \,  i\to \infty. \]
\end{proposition}
\begin{proof} The assumptions on $\widetilde F$ and $F_i$ imply that for
each $\eps>0$, each $k^*\in\N$ and each sufficiently large $i\in\N$, for
each trajectory $\bx_i$ of $F_i$ we can find a trajectory $\bx^\eps$
of $\widetilde F$ with $\|x^\eps_k-x_k\| \le \eps$,
$k=0,\ldots,k^*$. Hence, up to the time $k^*$ the trajectory
$\bx_\eps$ is also a trajectory of $F_\eps$. Thus, replacing $F$ by
$F_\eps$ we can follow the proof of Proposition \ref{prop:convergence}
in order to obtain 
\[ V_{(F_\eps,G)}(x) \le V_{(F_i,G_i)}(x) + 5\eps \]
for all sufficiently large $i\in\N$ and all $x\in K$. Now
(\ref{Vcont}) implies the assertion.
\end{proof}

\begin{remark} Basically, the continuity assumption (\ref{Vcont})
demands that an arbitrarily small relaxation of the state space
constraints does not lead to large changes in the optimal value
function. If $V_{(F,G)}$ is continuous on $K$ then one can expect  
(\ref{Vcont}) to hold for $p=\infty$ while if $V_{(F,G)}$ is
discontinuous on $K$ (note that state space restrictions may introduce
discontinuities in the optimal value function) then we would only
expect (\ref{Vcont}) to hold with $p<\infty$ because the location of the
discontinuity is likely to change when the state constraint changes. 
We conjecture that (\ref{Vcont}) holds under mild regularity
conditions on the optimal control problem, a formal verification,
however, is beyond the scope of this paper. 

In any case, we would like to emphasize that our result allows for a
rigorous convergence proof of the approximating multivalued game in
the presence of discontinuities, a feature which is rarely found in
other approximation techniques.
\end{remark}

\section{Discretization of the game}
\label{sec:discretization}

In this section we describe the set oriented discretization technique
which transforms our problem into a graph theoretic problem. 
In order to introduce our method, we first recall the corresponding
procedure for unperturbed systems developed in \cite{JungeOsinga:04}
before we turn to the general setting.

\subsection{Discretizing the Unperturbed System}

If $X$ is finite and there are no perturbations, then one can use a
shortest path algorithm like Dijkstra's method \cite{Dijkstra:59}, see also the appendix, in order to compute the value function, see, e.g., \cite{Bertsekas:95}.  In \cite{JungeOsinga:04}  it has been shown how to discretize general optimal control problems with continuous state space such that this approach can be applied.  We review this method here in a different formulation that directly carries over to the case of a perturbed control system in the next section. 

We consider a single valued control system $f:X\times U\to X$ ($f$
continuous, $X\subset\R^d$ and $U\subset\R^m$ compact, $0\in X$, $0\in U$,
$f(0,0)=0$), together with a continuous cost function $g:X\times U\to
[0,\infty)$ with $g(x,u) > 0$ for $x\neq 0$ and $g(0,0)=0$.  Let $\cP$
be a finite partition of $X$, i.e.\ $\cP$ is a finite set of mutually
disjoint subsets $P\subset X$.  Define the map $\pi:X\to \cP$, $\pi(x)=P$,
$x\in P$, as well as $\rho:X\rightrightarrows X$, $\rho=\pi^{-1}\circ\pi$ (i.e.\ to each $x$, $\rho$ associates the set of the partition $\cP$ which
contains $x$).

\paragraph{Box-enclosure of the system.}

Consider the multivalued game (which is actually a multivalued control system since there are no perturbations here) $(F,G)$ with
\[
F(x,u,w)=F(x,u):=\rho(f(x,u))\quad\text{and}\quad G(x,x_{1},u,w)=g(x,u).
\]
The optimality principle (\ref{multi-optprinc}) in this case reads
\begin{equation}\label{eq:optprinc unperturbed}
V_{(F,G)}(x) = \inf_{u\in U}\left\{g(x,u) + \inf_{x_{1}\in F(x,u)}V_{(F,G)}(x_{1})\right\}.
\end{equation}

\paragraph{Projection onto piecewise constant functions.} 

The right hand side of (\ref{eq:optprinc unperturbed}) defines an operator on real valued functions on $X$, the  \emph{dynamic programming operator} $L:\R^X\to\R^X$,
\[
L[v](x)=\inf_{u\in U}\left\{g(x,u)+\inf_{x_{1}\in F(x,u)} v(x_{1})\right\}.
\]
Note that the optimal value function $V_{(F,G)}$ is, by definition of $L$, a fixed point of $L$, i.e.\ $L[V_{(F,G)}]=V_{(F,G)}$.
Abusing notation, we identify the space $\R^\cP$ with the subspace of real valued functions on $X$ that are piecewise constant on the elements of the partition $\cP$ (in fact, we view $v\in\R^\cP$ as the function $v\circ\pi\in \R^X$).  We define the projection $\varphi:\R^X\to\R^\cP\subset\R^X$, 
\[
\varphi[v](x) = \inf_{x'\in\rho(x)} v(x'),
\]
and the corresponding \emph{discretized dynamic programming operator}
$L_{\cP}:\R^\cP\to\R^\cP$,
\[
L_{\cP} = \varphi\circ L.
\]
Explicitely, the discretized operator reads 
\begin{eqnarray*}
L_{\cP}[v](x) 
&=& \inf_{x'\in\rho(x)} \left\{\inf_{u\in U}\left\{g(x',u)+ \inf_{x_{1}\in F(x',u)} v(x_{1})\right\}\right\}\\
&=& \inf_{x'\in\rho(x),u\in U} \left\{g(x',u) + v(f(x',u))\right\},
\end{eqnarray*}
since $v\in\R^\cP$ is constant on each element of $\cP$, i.e.\ on each set $F(x',u)$.

We define the discretized optimal value function $V_{\cP}\in\R^\cP$ as
the unique fixed point of $L_{\cP}$ with 
$V_{\cP}(0)= 0$. Then 
$V_{\cP}$ satisfies the optimality principle 
\begin{equation}
\label{eq:optimality principle for the discretized system}
V_{\cP}(x) = \inf_{x'\in\rho(x),u\in U} \left\{g(x',u)+V_{\cP}(f(x',u))\right\}.
\end{equation}

\paragraph{Graph theoretic formulation.}

Note that since $\cP$ is finite, $V_{\cP}(f(x',u))$ in (\ref{eq:optimality principle for the discretized system}) can only take finitely many values. We can therefore rewrite (\ref{eq:optimality principle for the discretized system}) as
\begin{eqnarray}
V_{\cP}(x) &=& \min_{P\in \pi(f(\rho(x),U))}\;\inf_{x'\in\rho(x),u\in U: f(x',u)\in P} \left\{g(x',u)+V_{\cP}(P)\right\}
\label{eq:optprinc-enclosure}
\end{eqnarray}
where $V_{\cP}(P)=V_{\cP}(x)$ for any $x\in P\in\cP$.  If we define the multivalued map (or, equivalently, the directed graph)
$\cF:\cP\rightrightarrows\cP$, 
\begin{equation}
\cF(P) = \pi(f(\pi^{-1}(P),U)), \quad P\in\cP,
\end{equation}
and the cost function
\begin{equation}
\cG(P',P)=\inf\{g(x,u)\mid x\in P',f(x,u)\in P, u\in U\},
\end{equation}
we can rewrite (\ref{eq:optprinc-enclosure}) as
\[
V_{\cP}(P)=\min_{P_{1}\in\cF(P)}\{ \cG(P,P_{1}) + V_{\cP}(P_{1})\}.
\]
Note that this optimality principle can be interpreted as being solved by Dijkstra's algorithm.

\subsection{Discretization of the Perturbed System}

Now we want to carry over the discretization procedure from the last
section to our game setting. We proceed in a completely analogous way, additionally incorporating the perturbations now.  This will ultimately lead to a directed hypergraph (actually a \emph{forward hypergraph} or \emph{$F$-graph} in the terminology of \cite{gallo92directed}) instead of an ordinary graph for which we formulate the associated shortest path algorithm at the end of the section.

\paragraph{Box-enclosure of the system.} 

Consider the multivalued game $(F,G)$ with
\begin{equation}\label{dismultigame}
F(x,u,w) = \rho(f(x,u,w)) 
\quad\mbox{and}\quad
G(x,x_{1},u,w) = g(x,u),
\end{equation}
(where $f$ and $g$ are the control system and cost function introduced in Section~\ref{sec:problem}).  From the optimality principle~(\ref{multi-optprinc}) we obtain
\begin{eqnarray*}
\nonumber V_{(F,G)}(x) &=& \inf_{u\in U} \sup_{w\in W} \inf_{x_{1}\in F(x,u,w)}\left\{ g(x,u) +  V_{(F,G)}(x_{1})\right\}\\
&=& \inf_{u\in U} \left\{ g(x,u) +  \sup_{w\in W} \inf_{x_{1}\in F(x,u,w)} V_{(F,G)}(x_{1})\right\}.
\end{eqnarray*}

\paragraph{Projection onto piecewise constant functions.} 

The dynamic programming operator $L:\R^X\to\R^X$ here reads
\[
L[v](x)=\inf_{u\in U}\left\{g(x,u)+\sup_{w\in W} \inf_{x_{1}\in F(x,u,w)} v(x_{1})\right\}.
\]
Correspondingly, the discretized operator $L_\cP:\R^\cP\to\R^\cP$ is given by 
\begin{eqnarray*}
L_{\cP}[v](x) 
&=& \inf_{x'\in\rho(x)} \left\{\inf_{u\in U}\left\{g(x',u)+\sup_{w\in W} \inf_{x_{1}\in F(x',u,w)} v(x_{1})\right\}\right\}\\
&=& \inf_{x'\in\rho(x),u\in U} \left\{g(x',u)+\sup_{x_{1}\in F(x',u,W)} v(x_{1})\right\},
\end{eqnarray*}
since $v\in\R^\cP$ is constant on each element of $\cP$, i.e.\ on each set $F(x',u,w)$.

We define the discretized optimal value function $V_{\cP}\in\R^\cP$ as
the unique fixed point of $L_{\cP}$ with $V_{\cP}(P)= 0$ for all
partition elements $P\in\cP$ with 
$\overline{\pi^{-1}(P)}\cap O \ne \emptyset$. 
Then $V_{\cP}$ satisfies the optimality principle 
\begin{equation}
\label{eq:optimality principle for perturbed discretized system}
V_{\cP}(x) = \inf_{x'\in\rho(x),u\in U} \left\{g(x',u)+\sup_{x_{1}\in F(x',u,W)} V_{\cP}(x_{1})\right\}.
\end{equation}

\paragraph{Graph theoretic formulation.}

\begin{figure}[th]
\begin{center}
\includegraphics[width=0.9\textwidth]{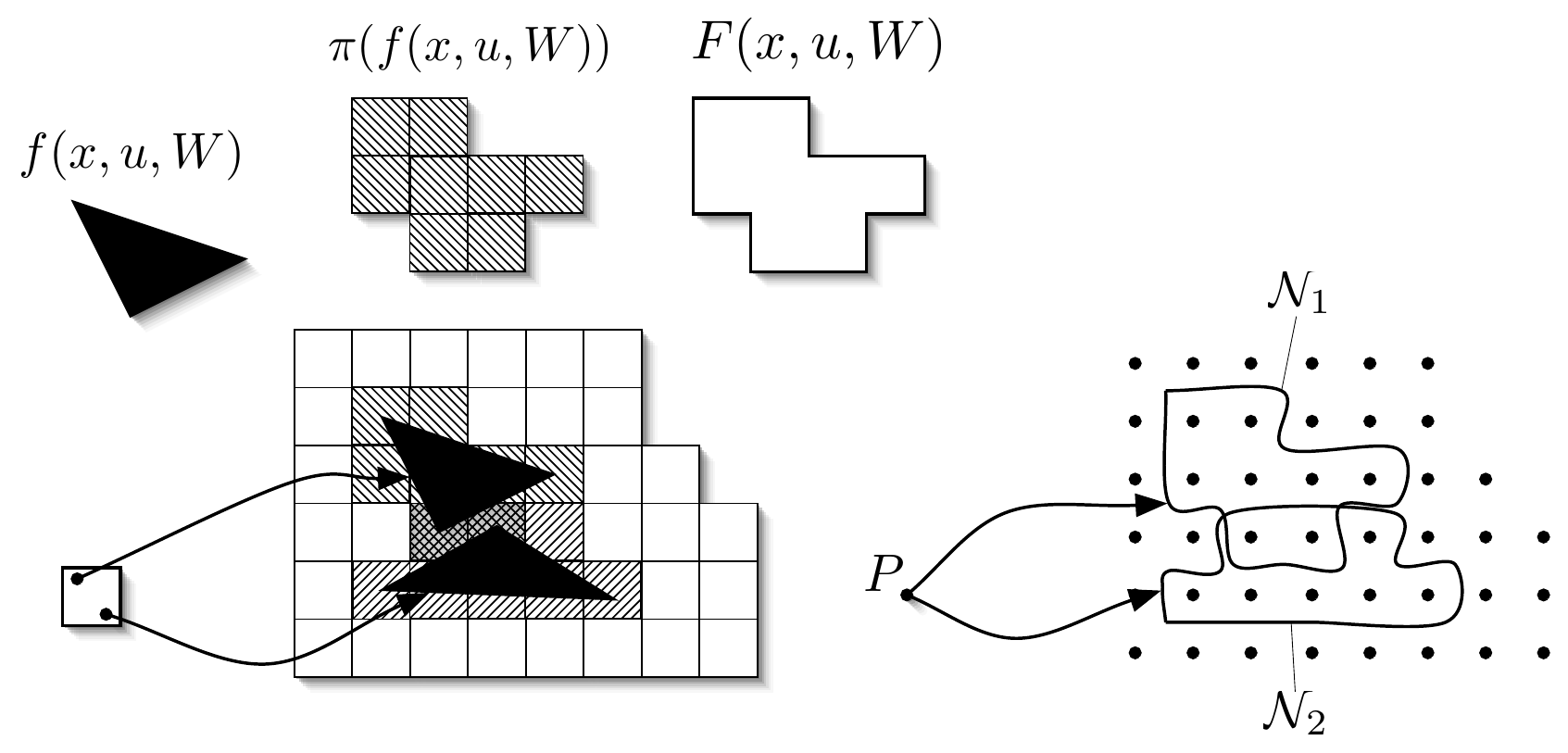}
\caption{Illustration of the construction of the hypergraph.}
\label{fig:hypergraph}
\end{center}
\end{figure}

In order to derive the corresponding shortest path algorithm, it is useful to formulate~(\ref{eq:optimality principle for perturbed discretized system}) equivalently in terms of an associated graph.  To this end note that for any pair $(x,u)\in X\times U$, the set $F(x,u,W)\subset X$ is the union of a finite set of elements from the partition $\cP$.  In particular, the family $\{ F(x',u,W) : (x',u)\in \rho(x)\times U\}$ of subsets of $X$ is finite for any $x\in X$.  Putting this in terms of a corresponding map on $\cP$: each partition element $P$ is mapped to a finite family $\{\cN_{i}\}_{i=1,\ldots,i(P)}$, $\cN_{i}\subset\cP$, of subsets of $\cP$ under all perturbations.  Formally, we have a directed hypergraph $(\cP,E)$ with the set $E\subset \cP\times 2^\cP$ of hyperedges given by
\[
E=\left\{(P,\cN) \mid \pi(F(x,u,W))=\cN \mbox{ for some } (x,u)\in P\times U\right\},
\]
or, equivalently, the multivalued map $\cF:\cP\rightrightarrows 2^\cP$,
\[
\cF(P) = \{ \pi(F(x,u,W)) : (x,u)\in P\times U\},
\]
c.f.\ Figure~\ref{fig:hypergraph}. 

If we define weights on the edges of this hypergraph by 
\[
\cG(P,\cN) = \inf \{ g(x,u) : (x,u)\in P\times U, \pi(F(x,u,W))=\cN\},
\]
then we can write (\ref{eq:optimality principle for perturbed discretized system}) equivalently as
\begin{equation}
\label{eq:perturbed graph OP}
V_{\cP}(P) = \inf_{\cN\in \cF(P)}\left\{ \cG(P,\cN) + \sup_{N\in\cN} V_{\cP}(N)\right\}.
\end{equation}

\paragraph{Dijkstra's method for the perturbed system.}

We are now going to generalize Dijkstra's algorithm (see the appendix) such that it computes the value function of a weighted directed hypergraph (i.e.\ the function defined by the optimality principle (\ref{eq:perturbed graph OP})).

Let $(\cP,E)$, $E\subset \cP\times 2^\cP$, be a hypergraph with weights $\cG:E\to [0,\infty)$.  In order to adapt Algorithm~\ref{dijkstra}, we need to modify the relaxing step in lines 7--9, such that the maximization over all perturbations (i.e.\ over $N\in\cN$) in (\ref{eq:perturbed graph OP}) is taken into account.   The modified version of lines 7--9 reads:
\begin{tabbing}
100\=mm\=mm\=mm\=mm\=  \kill
7\>\> for each $(Q,\cN)\in E$ with $P\in\cN$ \\
8\>\> \> if $V(Q) > \cG(Q,\cN) + \max_{N\in \cN} V(N)$ then \\
9\>\> \> \> $V(Q) := \cG(Q,\cN) + \max_{N\in \cN} V(N)$ 
\end{tabbing}

As justified by Proposition~\ref{dijkstra_max} (see the Appendix), if
$\cN\subset \cP\backslash \cQ$, then  
\[
\max_{N\in \cN} V(N) = V(P),
\]
and the node $Q$ will never be relaxed again.  On the other hand, if $\cN\not\subset \cP\backslash \cQ$, then $Q$ will be relaxed at a later time again  and we do not need to relax it in this iteration of the while-loop.  These considerations lead to the following further modification of lines 7--9:
\begin{tabbing}
100\=mm\=mm\=mm\=mm\=  \kill
7\>\> for each $(Q,\cN)\in E$ with $P\in\cN$\\
8\>\> \> if $\cN\subset \cP\backslash \cQ$ then \\
9\>\> \> \>if $V(Q) > \cG(Q,\cN) + V(P)$ then \\
10\>\> \> \> \> $V(Q) := \cG(Q,\cN) + V(P)$     
\end{tabbing}

Including the adapted initialization, the overall algorithm for the case of a perturbed system reads as follows.  Here, $\cD\subset\cP$ is the set of destination nodes which typically will be chosen as $\cD=\{P\in\cP : P\cap O\neq\emptyset\}$ (with the robust forward invariant set $O$ from Section~\ref{sec:problem}).

\begin{algorithm}\label{perturbed Dijkstra}
\textsc{Perturbed Dijkstra}$((\cP,E),\cG,\cD)$
\begin{tabbing}
100\=mm\=mm\=mm\=mm\=  \kill
1\>for each $P\in \cP$ set $V(P):=\infty$ \\
2\>for each $P\in \cD$ set $V(P) := 0$ \\
3\>$\cQ := \cP$\\
4\>while $\cQ\neq\emptyset$\\
5\>\> $P := \argmin_{P'\in \cQ} V(P')$\\
6\>\> $\cQ := \cQ\backslash\{P\}$\\
7\>\> for each $(Q,\cN)\in E$ with $P\in\cN$\\
8\>\> \> if $\cN\subset \cP\backslash \cQ$ then \\
%LG Tippfehler in Zeile 8 ausgebessert 
%(dort stand $\cN\subset \cP\backslash Q$)
9\>\> \> \>if $V(Q) > \cG(Q,\cN) + V(P)$ then \\
10\>\> \> \> \> $V(Q) := \cG(Q,\cN) + V(P)$     
\end{tabbing}
 \end{algorithm}
 
 We note that this algorithm bears similarities with the SBT-algorithm in \cite{gallo92directed}.  However, in our case the graph has a special structure (namely, the heads of the hyperedges consist of only a single node, i.e. we have an $F$-graph as defined in \cite{gallo92directed}).   This yields the subquadratic complexity in the number of nodes as derived above and thus gives an improvement over SBT.

\section{Implementation and Numerical Examples}\label{sec:implementation}

\subsection{Implementation}\label{ssec:implementation}

In the numerical realization we always let the state space $X$ be a box in $\R^d$ and
construct a partition $\cP$ of it by dividing $X$ uniformly into smaller
boxes.  In fact,  we realize this division by repeatedly bisecting the 
current division (changing the coordinate direction after each bisection).
The resulting sequence of partitions can efficiently be stored
as a binary tree --- see \cite{DellnitzHohmann:97} for more details.

In order to compute (or rather approximate) the set $E\subset \cP\times 2^\cP$ of hyperedges, we choose  finite sets $\tilde P\subset P$, $\tilde U\subset U$ and $\tilde W\subset W$ of \emph{test points} -- typically on an equidistant grid in each of these sets.  We then compute
\[
\tilde \cF(P) := \{\pi(F(x,u,\tilde W)) : (x,u)\in \tilde P\times\tilde U\} \subset 2^\cP
\]
as an approximation to $\cF(P)$ and correspondingly approximate the weights on the hyperedges by
\[
\tilde\cG(P,\cN) = \min \{ g(x,u) : (x,u)\in \tilde P\times \tilde U, \pi(F(x,u,\tilde W))=\cN\}.
\]

\paragraph{Time and space complexity.}  The time complexity of the standard Dijkstra algorithm (Algorithm~\ref{dijkstra} in the appendix) strongly depends on the data structure which is used in order to store the set $\cQ$. In particular, the complexity of the operations in lines 5 (extracting the node with minimal $V$-value) and line 9 (decreasing the $V$-value and the associated reorganization of the data structure) have a crucial influence.  In our implementation we are using a binary heap in order to store $\cQ$ which leads to a complexity of $\cO((|\cP|+|E|)\log|\cP|)$.

In the perturbed case (Algorithm~\ref{perturbed Dijkstra}), each hyperedge is considered at most $N$ times in line 7, with $N$ being a bound on the cardinality of the hypernodes $\cN$.  Additionally, we need to perform the check in line 8, which has linear complexity in $N$. Thus, the overall complexity of the perturbed Dijkstra algorithm is $\cO(|\cP|\log|\cP| + |E|N(N+\log|\cP|))$.

The space requirements grow linearly with the number of partition elements. Since typically the whole state space has to be covered, this number grows exponentially with the dimension of phase space (assuming a uniform partioning).  The concrete storage consumption strongly depends on the properties of the underlying control system.  While the number of hyperedges is essentially determined by the Lipschitz constant of $f$, the size of the hypernodes $\cN$ will crucially be influenced by the size of the perturbation.  In the applications that we have in mind in this paper, these numbers are of moderate size.

As a rule of thumb, the main computational effort in our approach goes into the construction of the hypergraph via the mapping of test points -- in particular, if the system is given by a short-time integration of a continuous time system.  Note that this ``sampling'' of the system will be required in any method that computes the value function.  Typically however, in standard methods like value iteration, certain points are sampled multiple times which leads to a higher computational effort in comparison to our approach.

\subsection{Numerical Examples}

\begin{figure}
	\centering
		\includegraphics[width=0.7\textwidth]{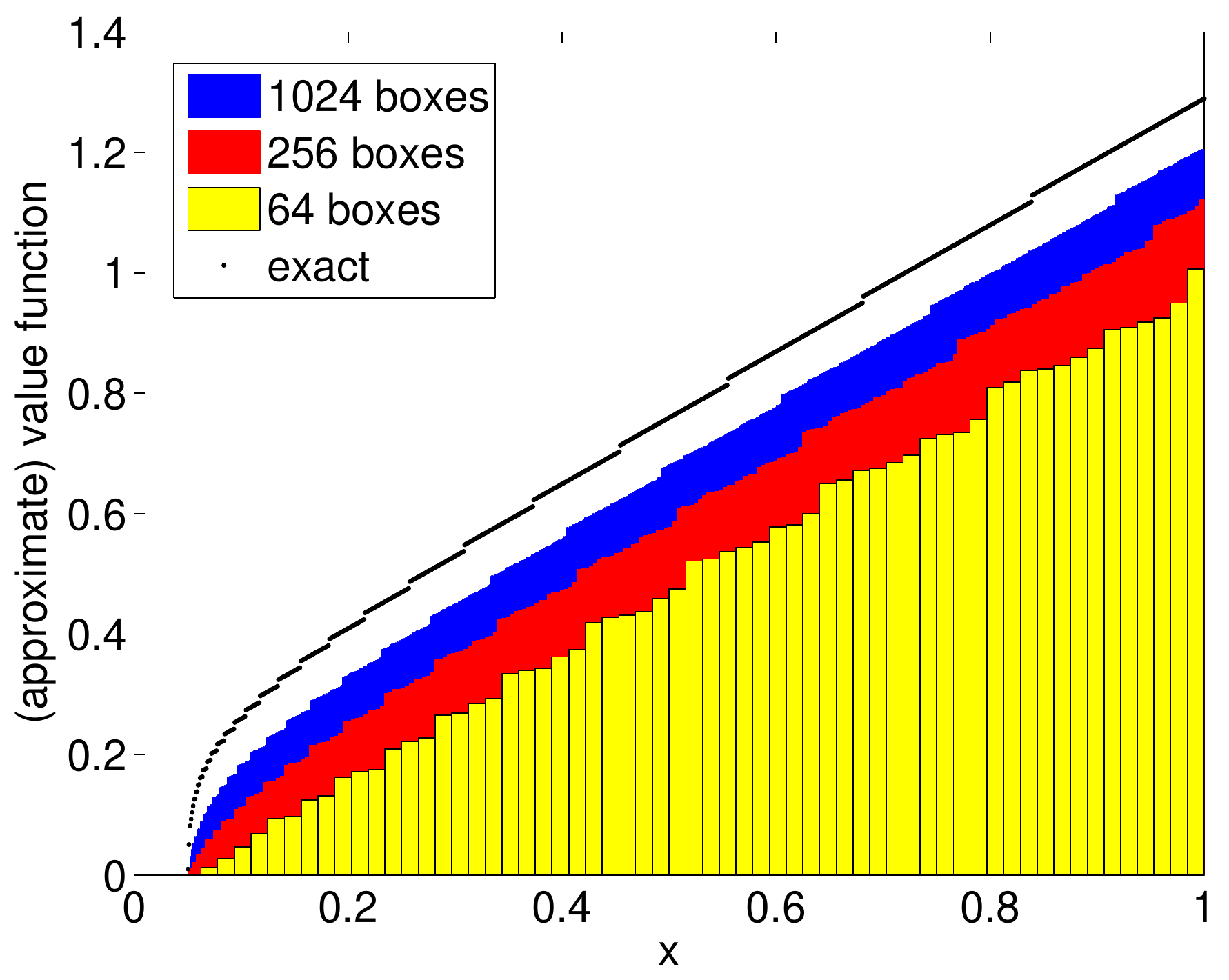}
		\caption{Perturbed simple 1D map: Upper value function and its approximations on various partitions.}
	\label{fig:v_simple}
\end{figure}

\paragraph{A simple 1D system.}  We start by looking at an additively perturbed version of a simple 1D map from \cite{GrJu:05}:
\[
x_{k+1} = x_k + (1-a)u_k x_k + w_k, \quad k=0,1,\ldots,
\]
with $x_k\in [0,1]$, $u_k\in [-1,1]$, $w_k\in [-\eps,\eps]$ for some $\eps > 0$ and the fixed parameter $a\in (0,1)$.  The cost function is
\[
g(x,u) = (1-a)x
\]
so that (regardless of how the perturbation sequence is chosen) the optimal control policy is to steer to the origin as fast as possible, i.e. to choose $u_k=-1$ for all $k$.  Similarly, the optimal strategy for the ``perturbing player'' is to slow down the dynamics as much as possible, corresponding to $w_k=\eps$ for all $k$.  The resulting dynamical system is the affine linear map
\[
x_{k+1} = ax_k + \eps, \quad k=0,1,\ldots,
\] 
which has a fixed point at $x=\eps/(1-a)$, i.e.\ under worst case conditions (assuming $w_k=\eps$ for all $k$) it will be impossible to get any closer than $\alpha_0:=\eps/(1-a)$ to the origin.  Correspondingly, we choose a neighborhood $O=[0,\alpha]$ with $\alpha > \alpha_0$ as our target region.  With 
\[
k(x) = \left\lceil \frac{\log\frac{\alpha-\alpha_0}{x-\alpha_0}}{\log a}\right\rceil + 1,
\]
the exact optimal value function is
\[
V(x) = (x-\alpha_0)\left(1-a^{k(x)}\right)+\eps k(x),
\]
as shown in Figure~\ref{fig:v_simple} for $a=0.8$, $\eps = 0.01$ and
$\alpha=1.1\alpha_0$.  In that Figure, we also show the approximate
optimal value functions on partitions of $64, 256$ and $1024$
intervals, respectively.  In the construction of the hypergraph, we
used an equidistant grid of ten points in each partition interval, in
the control space and in the perturbation space.  

%In our current (suboptimal) implementation, the computation on the finest partition took a couple of minutes on a recent workstation.

%LG: auskommentiert, weil dies zu pessimistisch klingt 
%    (die C-Version sollte doch viel schneller sein)

\paragraph{The inverted pendulum -- reloaded.} 
%LG: Hier habe ich einiges umformuliert und ergänzt; schau doch mal ob
%    das für Dich so ok ist
As a more challenging test case, we reconsider the problem of designing an optimal globally stabilizing controller for an inverted pendulum on a cart (see \cite{JungeOsinga:04, GrJu:05}): 
\begin{equation}\label{eq:pendulum}
        \left(\frac{4}{3} - m_r \cos^2 \varphi\right) \ddot\varphi + \frac{1}{2} m_r \dot\varphi^2 \sin 2 \varphi - \frac{g}{\ell}\sin \varphi  = - u\;\frac{m_r}{m
        \ell} \cos \varphi.
\end{equation}
The equation models the (planar) motion of an inverted pendulum with mass $m=2$
on a cart with mass $M=8$ which moves under an applied horizontal force $u$.
The angle $\varphi$ measures the offset angle from the vertical up position. 
The parameter $m_r = m / (m + M)$ is the mass ratio and $\ell = 0.5$ the distance of the pendulum mass from the pivot. We use $g = 9.8$ for the gravitational constant.  The instantaneous cost is 
\begin{equation}
\label{eq:costcont}
        q(\varphi,\dot\varphi,u) = \frac{1}{2} \left(0.1 \varphi^2 + 0.05 \dot\varphi^2 + 0.01 u^2\right).
\end{equation}
Denoting the evolution operator of the control system
(\ref{eq:pendulum}) for constant control functions $u$ by
$\Phi^t(t,u)$, we consider the time-$T$-map $\Phi^T(x,u)$ of this
system as our 
discrete time system with $T=0.1$. The map $\Phi^T$ is approximated
via the classical Runge-Kutta scheme of order $4$ with step size
$0.02$. Thus we arrive at the cost function   
\[
  g(\varphi,\dot\varphi,u) = \int_0^T q(\Phi^t((\varphi,\dot\varphi),u),u)\; dt,
\]
We choose $X = [-8, 8] \times [-10,10]$ as the region of interest.

\begin{figure}
	\centering
		\includegraphics[width=0.49\textwidth]{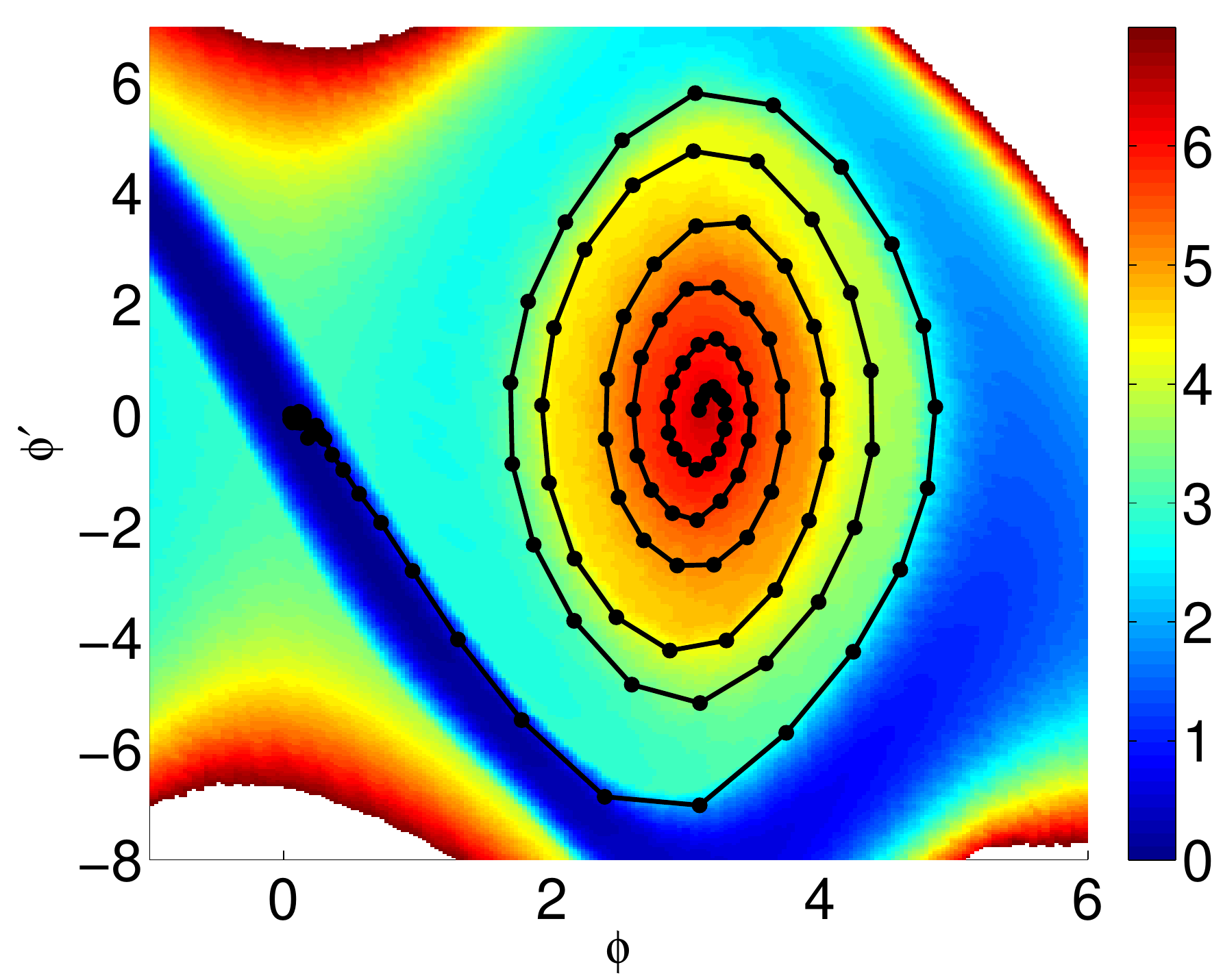}\hfill
		\includegraphics[width=0.47\textwidth]{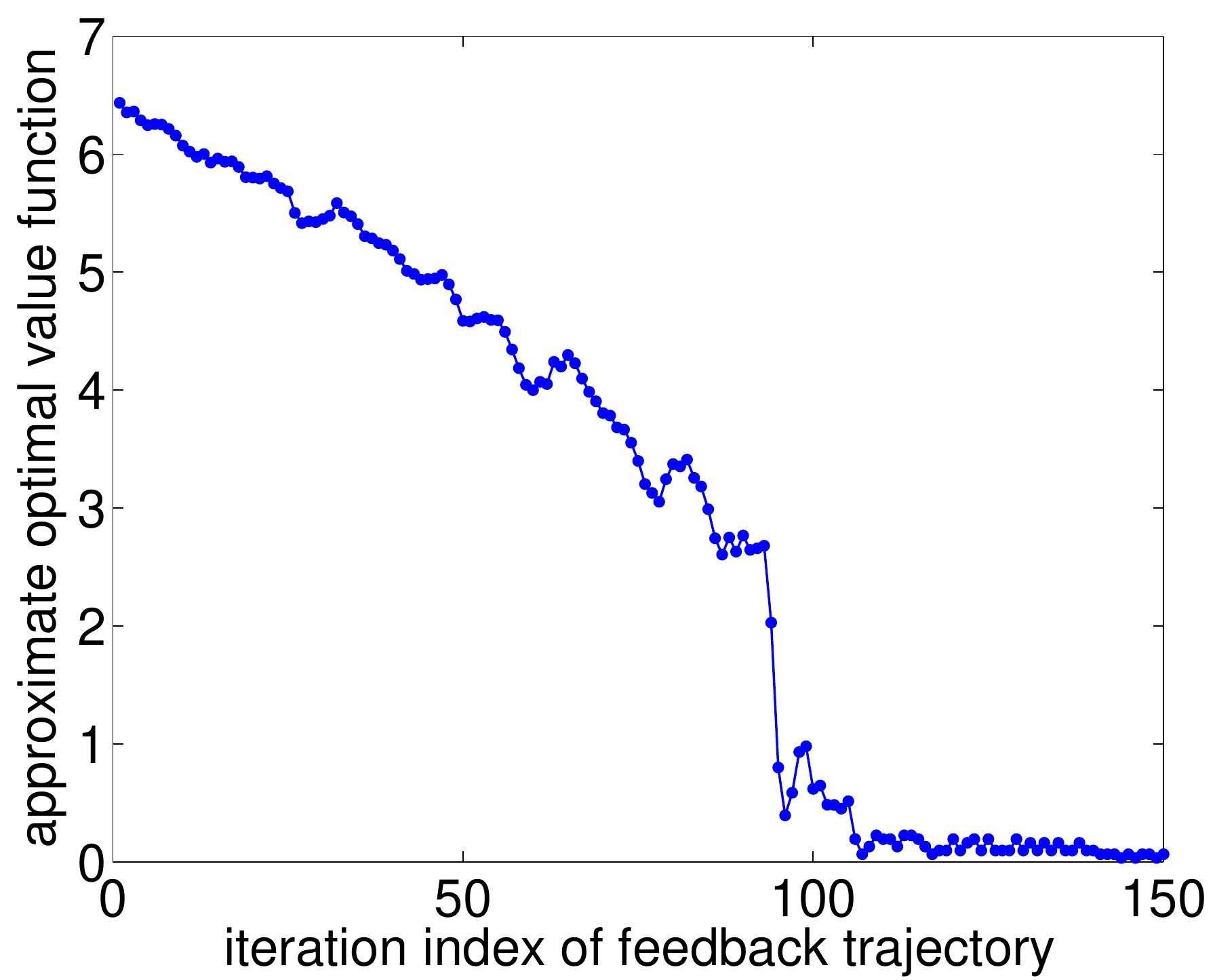}
		\caption{Approximate optimal value function and feedback trajectory (left) and the approximate optimal value function along the feedback trajectory (right) for the inverted pendulum on a $2^{18}$ box partition.}
	\label{fig:v_decrease}
\end{figure}

In \cite{GrJu:05}, a feedback trajectory with initial value $(3.1,0.1)$
was computed that was based 
on an approximate optimal value function on a partition of $2^{18}$
boxes (cf.\ Figure~\ref{fig:v_decrease} (left)).  In contrast to what
one might expect, the approximate optimal value function does actually
not decrease monotonically along this trajectory (cf.\
Figure~\ref{fig:v_decrease} (right)). 
This effect is due to the fact that the discretization method used in
\cite{GrJu:05} allows for jumps in the trajectories which cannot be
reproduced by the real system. The fact that the approximate optimal
value function is not always decreasing indicates that the
approximation accuracy in this example is just fine enough to allow
for stabilization, and in fact, on a coarser partition of $2^{14}$
boxes, the associated feedback is not stabilizing this initial
condition any more.   

We are now going to use the approach developed in this paper in order
to design a stabilizing feedback controller on basis of the coarser
partition ($2^{14}$ boxes).  To this end, we imagine the perturbation
of our system being given as ``for a given state
$(\varphi,\dot\varphi)$, be prepared to start anywhere in the box that
contains $(\varphi,\dot\varphi)$'', i.e.\ we define our game by 
\[
F((\varphi,\dot\varphi),u,W) := \Phi^T(B,u),
\] 
where $B\in\cP$ is the box in the partition $\cP$ under consideration which contains the point $(\varphi,\dot\varphi)$.  Note that we do not need to parameterize the points in $\Phi^T(B,u)$ with $w\in W$ for the construction of the hypergraph.   

\begin{figure}
	\centering
		\includegraphics[width=0.49\textwidth]{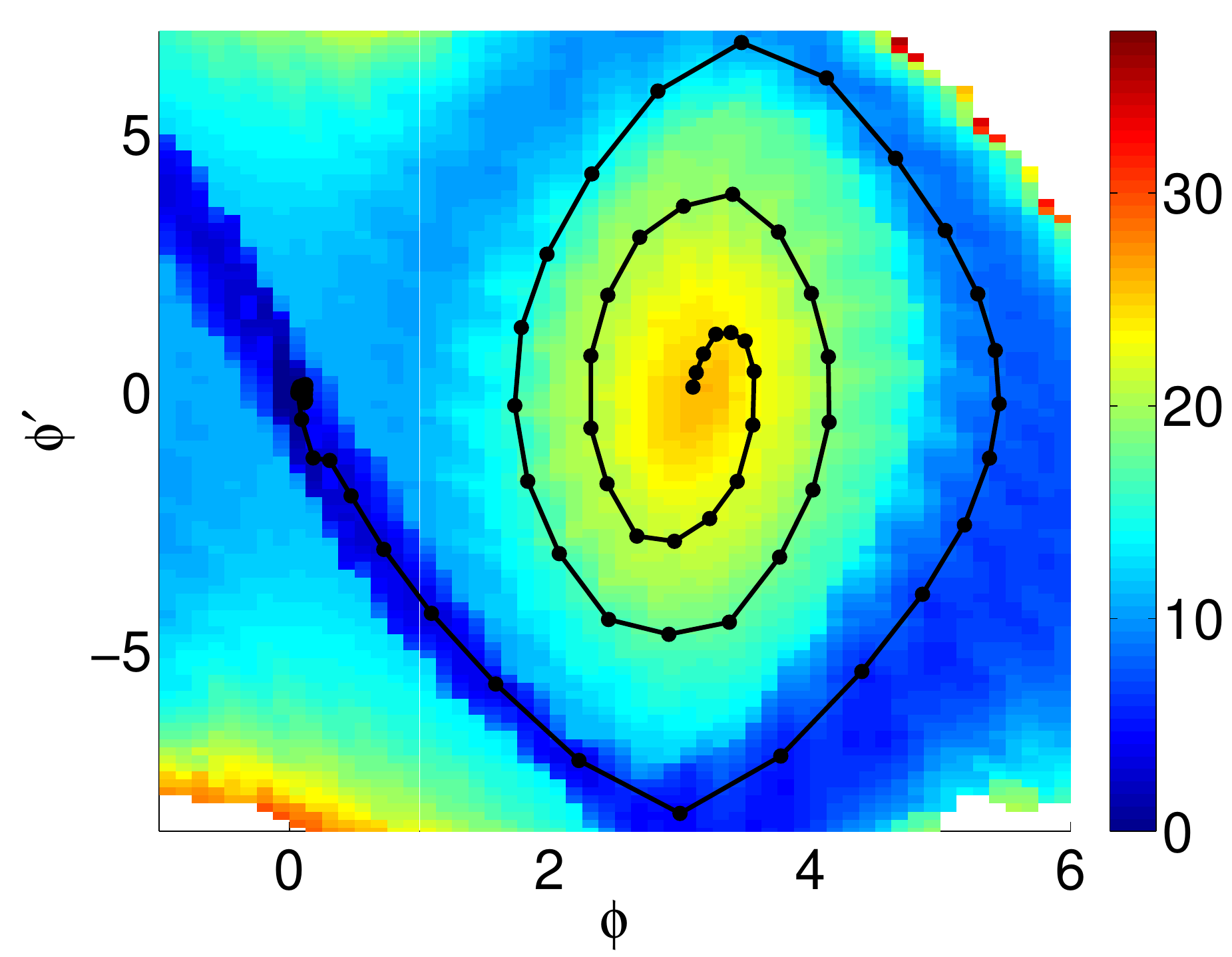}\hfill
		\includegraphics[width=0.47\textwidth]{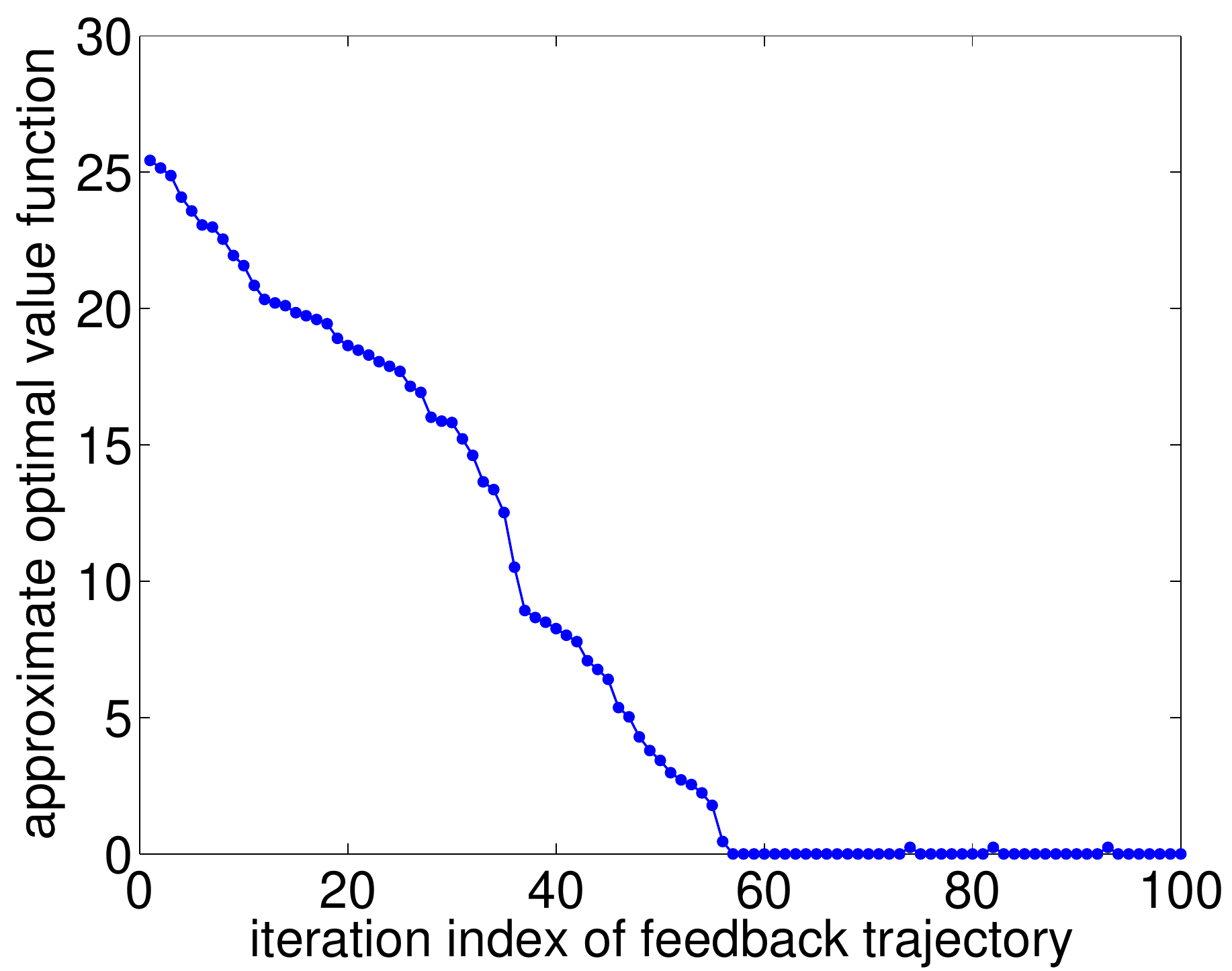}
		\caption{Approximate upper value function and feedback trajectory (left) and the approximate upper value function along the feedback trajectory (right) for the inverted pendulum on a $2^{14}$ box partition using the robust feedback construction.}
	\label{fig:v-robust}
\end{figure}

Figure~\ref{fig:v-robust} shows the approximate upper value function
on a partition of $2^{14}$ boxes with target region $O=[-0.1,0.1]^2$ as well as the trajectory generated
by the associated feedback for the initial value $(3.1,0.1)$.   As
expected, the approximate value function is decreasing monotonically
along this trajectory. Furthermore, despite the fact that we used
considerably fewer boxes as for Figure~\ref{fig:v_decrease}, the
resulting trajectory is obviously closer to the optimal one because it
converges to the origin much faster. 

\section{Convergence Analysis}\label{sec:convergence}

In this section we show that and in which sense the approximate
optimal value function constructed in the preceeding section converges
to the true one as the underlying partitions are refined, using the
abstract results for multivalued games developed in the Sections
\ref{sec:mvg} and \ref{stateconstraints}.

We begin with the following observation on the relation between $V_\cP$ and $V_{(F,G)}$ with $F$, $G$ from (\ref{dismultigame}).

\begin{proposition}\label{prop:projection} Consider the discretized optimal value function
  $V_\cP$ and the optimal value function $V_{(F,G)}$ from
  (\ref{multival}) corresponding to the game
  (\ref{dismultigame}).  If $V_{(F,G)}$ is continuous on $\partial O$, then these functions are related by 
  \[ V_\cP(x) = \inf_{x'\in \rho(x)} V_{(F,G)}(x'). \] 
\end{proposition}
{\bf Proof:} First note that both functions are nonnegative. 
>From the previous considerations it follows that the
functions satisfy the optimality principles
\begin{equation}\label{opt1} 
V_{(F,G)}(x) = \inf_{u\in U} \sup_{w\in W} \inf_{x_{1}\in F(x,u,w)}\left\{
  g(x,u) +  V_{(F,G)}(x_{1})\right\} \end{equation}
and
\begin{equation}\label{opt2} 
V_\cP(x) = \inf_{x'\in \rho(x)} \inf_{u\in U} \sup_{w\in W}
\inf_{x_{1}\in F(x',u,w)}\left\{g(x',u) +
  V_\cP(x_{1})\right\}. \end{equation} 
In order to show 
\begin{equation}\label{desineq} 
\inf_{x'\in \rho(x)} V_{(F,G)}(x')\le V_\cP(x), \end{equation}
we number the elements $P_i$ of $\cP$ such that $i_2>i_1$ implies 
$V_\cP|_{P_{i_2}} \ge V_\cP|_{P_{i_1}}$. We first consider those
elements $P_i$, $i=1,\ldots, j$, for which we have $V_\cP|_{P_i}\equiv
0$ which by our assumptions on $V_\cP$ and $g(x,u)$ is equivalent to   
$\overline{\pi^{-1}(P_i)} \cap O\ne \emptyset$. 

In case that $\pi^{-1}(P_i) \cap O\ne \emptyset$, we can find 
$x_0\in \pi^{-1}(P_i)\cap O$ and $u_0\in U$ such that 
$F(x_0,u_0,w)\subset O$ for all $w\in W$. In particular,
for any fixed $w$ we
find $x_1\in F(x_0,u_0,w)\cap O$ for which we proceed the same way,
which yields $F(x_1,u_1,w)\subset O$ for all $w\in W$. Hence, given a
perturbation strategy $\beta(\bu)$ we find a control sequence $\bu$
such that $\mathcal{X}_{F}(x_0,\bu,\beta(\bu))\subset O$ implying 
\[ J_{(F,G)}(x_{0},\bu,\beta(\bu)) =
\inf_{(x_{k})_{k}\in\mathcal{X}_{F}(x,\bu,\beta(\bu))}\sum_{k=0}^\infty
G(x_{k},x_{k+1},u_{k},\beta(\bu)_{k}) = 0 \]
and thus 
\[ \inf_{x'\in \rho(x_{0})} V_{(F,G)}(x') \le V_{(F,G)}(x_0) = 0 \le 
V_\cP(x_{0}), \]
which shows (\ref{desineq}) for $\rho(x)=P_i$ with $\pi^{-1}(P_i) \cap O\ne \emptyset$.  In fact, what we showed is that $V_{(F,G)}(x)=0$ for $x\in O$.  Since we assumed that $V_{(F,G)}$ is continuous on $\partial O$, we also get
\[
\inf_{x'\in P_{i}} V_{(F,G)}(x')=0
\]
for $P_{i}$ with $\overline{\pi^{-1}(P_i)} \cap O\ne \emptyset$, but
$\pi^{-1}(P_i) \cap O = \emptyset$.

Now we proceed by induction over $i\ge j+1$. We pick some $i\ge j+1$
and assume that the desired inequality (\ref{desineq}) holds for 
$\rho(x)=P_1,\ldots, P_{i-1}$. We fix 
$x\in X$ with $\rho(x)=P_i$ and an arbitrary $\eps>0$. Then we pick
$x''\in P_i$ such that the infimum over $x'$ in (\ref{opt2}) 
is attained up to $\eps$. Thus we obtain
\begin{eqnarray*}
V_\cP(x) & = & \inf_{x'\in \rho(x)} \inf_{u\in U} \sup_{w\in W}
\inf_{x_{1}\in F(x',u,w)}\left\{g(x',u) + V_\cP(x_{1})\right\}\\
& \ge & \inf_{u\in U} \sup_{w\in W}
\inf_{x_{1}\in F(x'',u,w)}\left\{g(x'',u) + V_\cP(x_{1})\right\} -
\eps \\
& = & \inf_{u\in U} \sup_{w\in W}
\inf_{x_{1}\in F(x'',u,w)}\left\{g(x'',u) + V_{(F,G)}(x_{1})\right\} -
\eps \\
& = & V_{(F,G)}(x'') - \eps \, \ge \, \inf_{x'\in P_i} V_{(F,G)}(x') - \eps, 
\end{eqnarray*}
where we have used the induction assumption in the third step as
follows: the inequality $g(x,u)> 0$ implies 
$V_\cP(x_{1}) < V_\cP(x) = V_\cP|_{P_i}$, furthermore we have 
$x_1 \in F(x'',u,w) = P_{i'}$ for 
some $i'\in\N$, i.e., $V_\cP(x_{1}) = V_\cP|_{P_{i'}}$. This implies  
$V_\cP|_{P_i} > V_\cP|_{P_{i'}}$ and consequently $i>i'$. Hence by the
induction assumption we have 
\[ \inf_{x_{1}\in F(x'',u,w)}V_\cP(x_{1}) = V_\cP|_{P_{i'}} =
\inf_{x_{1}\in F(x'',u,w)}V_{(F,G)}(x_{1}). \]
Now, since $\eps>0$ was arbitrary, we obtain (\ref{desineq}).

The converse inequality $V_\cP(x)\le \inf_{x'\in\rho(x)} V_{(F,G)}(x)$ follows by a similar
induction argument using the fact that \eqref{opt1} always yields a
larger value than \eqref{opt2} due to the additional minimization over
$x'$ in \eqref{opt2}. \endproof

\begin{remark}
Note that in order to obtain the assertion from the preceeding proposition, it is sufficient that the union of those partition elements that have nonempty intersection with $O$ form a neighborhood of $O$.  If this is true, one can actually drop the assumption on the continuity of $V_{(F,G)}$ on $\partial O$.  
\end{remark}

We now consider a sequence of increasingly finer partitions of $X$ and ask under which conditions the corresponding approximate optimal value functions converge to the value function of the game $(f,g)$.  In a \emph{nested} sequence of partitions, each element of a partition is contained in an element of the preceding partition.

The following theorem states our main convergence result. It shows
that we obtain $L^\infty$ convergence on compact sets on which
$V_{(f,g)}$ is continuous and --- under a mild regularity condition on
the set of discontinuities --- $L^1$ convergence on every compact set
on which $V_{(f,g)}$ is bounded. We first consider problems without
state space constraints and address the constrained case in Remark
\ref{constrconv}, below.

\begin{theorem}\label{thm:convergence}
Let $(\cP_{i})_{i\in\N}$ be a nested sequence of partitions of $X$ such that
\[
\sup_{x\in X} H(\rho_{i}(x),\{x\})\to 0 \quad \text{as } i\to\infty.
\]
Assume that $g(x,u)$ is continuous, that $g(x,u)>0$ for 
$x\not\in O$ and that $V_{(f,g)}$ is
continuous on $\partial O$. Then 
\[
\|V_{\cP_{i}}|_{K_i}-V_{(f,g)}|_{K_i}\|_{\infty} \to 0 \quad \text{as }
i\to \infty 
\]
for every compact set $K\subseteq X$ on which $V_{(f,g)}$ is continuous and
$$K_i=\bigcup_{P\in\cP_i,\, \pi^{-1}(P)\subset K}\pi^{-1}(P)$$ being the
largest subset of $K$ which is a union of partition elements $P\in\cP_i$.

If we assume furthermore that the set of discontinuities of
$V_{(f,g)}$ has zero Lebesgue measure, then 
\[
\|V_{\cP_{i}}|_{K}-V_{(f,g)}|_{K}\|_{L^1} \to 0 \quad \text{as } i\to \infty
\]
on every compact set $K\subseteq X$ with 
$\sup_{x\in K}V_{(f,g)}(x)<\infty$. 
\end{theorem}

\begin{proof}
We use Proposition~\ref{prop:convergence} with $(F,G)=(f,g)$ ($f$
interpreted as a set valued map) and Proposition~\ref{prop:projection}. 

Note that since $F_{i}(x,u,w)=\rho_{i}(f(x,u,w))$ and
$G_{i}(x,u,w)=g(x,u)$, the games $(F_{i},G_{i})$ are enclosures of
$(f,g)$ (in fact, since the sequence of partitions is nested, for
every $i$, $(F_{i},G_{i})$ is an enclosure of $(F_{i+1},G_{i+1})$).
Under the assumptions of the theorem, all assumptions of
Proposition~\ref{prop:convergence} are satisfied.  In particular, by
the assumptions on $g$ and since $X$ and $U$ are compact, we know that
there exists a function $\alpha\in\KK_{\infty}$ such that 
\[
G_{i}(x,x_{1},u,w)=g(x,u)\geq \alpha(d(x,O) + d(x_{1},O))
\]
for all $i$.  Thus, $V_{(F_{i},G_{i})}$ converges uniformly to
$V_{(f,g)}$ on $K$.
In order to show the $L^\infty$ convergence on $K_i$ observe that if
$V_{(f,g)}$ is continuous on $K$ then it is also 
uniformly continuous on $K$ which implies 
\[ \sup_{P\in\cP_i,\, \pi^{-1}(P)\subset K}|\inf_{x\in P}
V_{(f,g)}(x) - \sup_{ x\in P} V_{(f,g)}(x)| \to 0\]
as $i\to\infty$. Thus we can use
Proposition~\ref{prop:projection} in order to conclude 
\begin{eqnarray*} 
\|V_{\cP_{i}}|_{K_i}-V_{(f,g)}|_{K_i}\|_{\infty} & \le &
\sup_{P\in\cP_i,\, \pi^{-1}(P)\subset K}
|V_{\cP_{i}}|_{P} - \sup_{x\in P} V_{(f,g)}(x) | \\
& = & \sup_{P\in\cP_i,\, \pi^{-1}(P)\subset K}
|\inf_{y\in P} V_{(F_i,G_i)}(y) - \sup_{x\in P} V_{(f,g)}(x) | \\
& \le & \sup_{P\in\cP_i,\, \pi^{-1}(P)\subset K}
\Big\{ |\inf_{y\in P} V_{(F_i,G_i)}(y) - \inf_{x\in P}
V_{(f,g)}(x) | \\
&& \hspace{2cm} 
+\,  |\inf_{x\in P} V_{(f,g)}(x) - \sup_{ x\in P} V_{(f,g)}(x)|\Big\}
\to 0 \end{eqnarray*}
as $i\to\infty$.

In order to show the $L^1$ convergence, observe that the uniform
convergence $V_{(F_{i},G_{i})}\to V_{(f,g)}$ on $K$ implies
\[
\|V_{(F_{i},G_{i})}|_{K}-V_{(f,g)}|_{K}\|_{L^1}\to 0\quad\text{as } i\to\infty.
\]
It thus remains to show that $V_{(F_{i},G_{i})}|_{K}-V_{\cP_{i}}|_{K}\to 0$ in $L^1$.  Let $D$ be the set of discontinuities of $V_{(f,g)}$ and $\cD_{i}=\{P\in\cP_{i},P\cap D\neq\emptyset\}$. We write
\begin{eqnarray*}
\int_{K} V_{(F_{i},G_{i})} - V_{\cP_{i}}\; dm& = & I_{i,1}+I_{i,2}
\end{eqnarray*}
with
\begin{eqnarray}
I_{i,1} &=& \sum_{P\in\cD_{i}} \int_{P\cap K} V_{(F_{i},G_{i})} - V_{\cP_{i}}\; dm ,\\ 
I_{i,2} &=& \sum_{P\in\cP_{i}\backslash \cD_{i}} \int_{P\cap K} V_{(F_{i},G_{i})} - V_{\cP_{i}}\; dm.
\end{eqnarray}
Because of $V_{(f,g)}\geq V_{(F_{i},G_{i})}$, the assumption that $D$ has zero Lebesgue measure and $H(\rho_{i}(x),\{x\})\to 0$, we have that $I_{i,1}\to 0$ for $i\to\infty$.  Using Proposition~\ref{prop:projection}, the compactness of $K$, and the fact that $V_{(F_{i},G_{i})}|_{K}\to V_{(f,g)}|_{K}$ uniformly, we also obtain that $I_{i,2}\to 0$ as $i\to\infty$, i.e.\ $V_{(F_{i},G_{i})}|_{K}-V_{\cP_{i}}|_{K}\to 0$ in $L^1$ and thus the assertion of the theorem.
\end{proof}

\begin{corollary}
Under the assumptions of Theorem~\ref{thm:convergence} we have
\[
V_{\cP_{i}}(x)\to V_{(f,g)}(x)\quad\text{ as } i\to\infty
\]
for Lebesgue-almost all $x\in K$, where $K$ is any compact subset of the domain of $V_{(f,g)}$.
\end{corollary}

\begin{proof}
By standard arguments, there exists a subsequence $(i(j))_{j}$ such that $V_{\cP_{i(j)}}(x)\to V_{(f,g)}(x)$ as $j\to\infty$ for Lebesgue-almost all $x\in K$.  Since $(V_{\cP_{i}}(x))_{i}$ is monotone, we obtain the assertion.
\end{proof}

\begin{remark}\label{constrconv} Using Proposition
\ref{prop:sc-convergence} instead of Proposition \ref{prop:convergence} 
it is easily seen that our convergence results remain valid 
in case of state space constraints if we assume condition
(\ref{Vcont}) for $\tilde F(x,u,w)=\{f(x,u,w)\}$.  
In this case, the first assertion of Theorem \ref{thm:convergence} 
will hold for the $p$--norm from (\ref{Vcont}) instead of the
$\infty$--norm. 
\end{remark}

\section{Feedback Construction}\label{sec:feedback}

As usual, we use the approximate optimal value function $V_{\cP}$ and the optimality principle (\ref{eq:optimality principle}) in order to construct an approximate optimal feedback.  More precisely,  for any point $x\in S_{0}$, $S_{0}:=\{ x\in X : V_{(f,g)}(x) < \infty\}$, we define
\[
u_{\cP}(x) = \argmin_{u\in U} \max_{w\in W} \;\{ g(x,u) + V_{\cP}(f(x,u,w)) \}.
\]
We can immediately adapt Theorem 3 from \cite{GrJu:05} in order to obtain a statement about the performance of this feedback.  The following result in particular shows that the feedback is robust with respect to arbitrary perturbations of the system.

\begin{theorem}
Let the assumptions of Theorem~\ref{thm:convergence} be satisfied. Let $D\subset S_{0}$ be an open set with compact closure, such that $\overline D\subset S_{0}$, $O\subset D$ and on which $V_{(f,g)}$ is continuous. Let $c>0$ be such that the inclusion $D_{c}(i_{0}):=V_{\cP_{i_{0}}}^{-1}([0,c])\subset D$ holds for some $i_{0}\in\N$.  Then there exists a function $\delta:\R\to\R$ with $\lim_{\alpha\to 0}\delta(\alpha)=0$ such that for all sufficiently small $\eps$, all sufficiently large $i$, all $\eta\in (0,1)$, all $x_{0}\in D_{c}(i)$ and all perturbation sequences $(w_{k})_{k}\in W^\N$, the trajectory generated by
\[
x_{k+1}=f(x_{k},u_{\cP_{i}}(x_{k}),w_{k})
\]
satisfies
\[
V(x_{k})\leq \max\left\{ V(x_{0}) - (1-\eta)\sum_{j=0}^{k-1} g(x_{j},u_{\cP_{i}}(x_{j})), \delta(\eps/\eta)+\eps\right\}.
\]
\end{theorem}

\begin{proof}
We only point out how to suitably modify the proof of Theorem 3 in \cite{GrJu:05}. First note that according to Theorem~\ref{thm:convergence}, $V_{\cP_{i}}$ converges uniformly to $V_{(f,g)}$ on $D$.  The second observation is that if we choose $i_{1}\in \N$, $i_{1}>i_{0}$ such that $V_{(f,g)}-V_{\cP_{i}}(x)\leq \eps/2$ for $i\geq i_{1}$ and all $x\in D_{c}(i_{1})$, then
\begin{eqnarray*}
V_{\cP_{i}}(x) + \eps/2  \;\geq\;   V(x) 
 & = & \inf_{u\in U}\sup_{w\in W}\{ g(x,u)+V(f(x,u,w))\}\\
 &\geq & \min_{u\in U}\max_{w\in W} \{ g(x,u)+V_{\cP_{i}}(f(x,u,w))\}\\
 &=&  g(x,u_{\cP_{i}}(x))+ \max_{w\in W}V_{\cP_{i}}(f(x,u_{\cP_{i}}(x),w)),
\end{eqnarray*}
i.e.
\[
V_{\cP_{i}}(x_{k+1}) \leq V_{\cP_{i}}(x_{k}) - g(x,u_{\cP_{i}}(x)) + \eps/2
\]
for all $x_{k+1}\in f(x_{k},u_{\cP_{i}}(x),W)$.  The rest of the proof of Theorem~3 in \cite{GrJu:05} remains the same.
\end{proof}

\begin{remark}
A particular application of our result is to robustify the feedback
construction from \cite{GrJu:05} with respect to small perturbations
which may be due, e.g., to discretization errors resulting from the
numerical computation of the discrete time system from an ordinary
differential equation. For this purpose, a particularly convenient way
is to consider an ``$\eps$-inflated'' system related to the original
unperturbed system.  More precisely, given an unperturbed
control system $f:X\times U\to X$, one considers the perturbed system 
\[
x_{k+1}=f(x_{k},u_{k})+\eps w_{k}, \quad k=0,1,\ldots,
\]
with $w_{k}\in [-1,1]^d$ for some (small) $\eps >0$.  In the numerical
realization, the sets $F(x,u,W)=f(x,u)+\eps [-1,1]^d$ are easy to
construct using ideas from rigorous discretization, see
\cite{Junge:00a,Gruene:02}.  
\end{remark}

\begin{appendix}

\section{Dijkstra's Method}

Let $(\cP,E)$ be a finite directed graph with edge weights $g:E\to [0,\infty)$.  Let $D\in \cP$ be the \emph{destination node}.  The following algorithm \cite{Dijkstra:59} computes the length $V(P)\in [0,\infty)$ of the shortest path from $P$ to $D$ for all nodes $P\in \cP$.

\begin{algorithm}\label{dijkstra}
\textsc{Dijkstra}$((\cP,E),g,D)$
\begin{tabbing}
100\=mm\=mm\=mm\=mm\=  \kill
1\>for each $P\in \cP$ set $V(P):=\infty$ \\
2\>$V(D) := 0$ \\
3\>$\cQ := \cP$\\
4\>while $\cQ\neq\emptyset$\\
5\>\> $P := \argmin_{P'\in \cQ} V(P')$\\
6\>\> $\cQ := \cQ\backslash\{P\}$\\
7\>\> for each $Q\in\cP$ with $(Q,P)\in E$ \\
8\>\> \> if $V(Q) > g(Q,P) + V(P)$ then \\
9\>\> \> \> $V(Q) := g(Q,P) + V(P)$     
\end{tabbing}
\end{algorithm}

An important feature of this algorithm is given by the following proposition, which follows immediately from the construction of the algorithm and the fact that the edge weights are nonnegative.

\begin{proposition}\label{dijkstra_max}
During the while-loop in lines 4-9 of Algorithm~\ref{dijkstra} it holds that
\[
V(P) \ge V(P')\quad\text{for all } P'\in \cP\backslash \cQ.
\]
\end{proposition}

\end{appendix}

\section*{Acknowledgement}

We thank Marcus von Lossow for helpful comments on the complexity analysis.

\baselineskip=12pt
\small
\bibliographystyle{unsrt}
\bibliography{GrJu04a}

\end{document}